\documentclass[11pt,english]{amsart}
\usepackage[T1]{fontenc}
\usepackage[latin1]{inputenc}
\usepackage{latexsym}
\usepackage{amsthm}
\usepackage{amssymb}

\providecommand{\tabularnewline}{\\}

\numberwithin{equation}{section} %% Comment out for sequentially-numbered
\numberwithin{figure}{section} %% Comment out for sequentially-numbered
\theoremstyle{plain}
\theoremstyle{plain}
\newtheorem{thm}{Theorem}
  \theoremstyle{plain}
  \newtheorem{prop}[thm]{Proposition}
  \theoremstyle{definition}
  \newtheorem{condition}[thm]{Hypothesis}
  \theoremstyle{plain}
  \newtheorem{lem}[thm]{Lemma}
  \theoremstyle{remark}
  \newtheorem{rem}[thm]{Remark}
  \theoremstyle{plain}
  \newtheorem{cor}[thm]{Corollary}
  \theoremstyle{definition}
  \newtheorem{defn}[thm]{Definition}
  \theoremstyle{remark}
  \newtheorem{notation}[thm]{Notation}

\AtBeginDocument{
  
}

\usepackage{babel}
\addto\extrasfrench{\providecommand{\fg}{\ifdim\lastskip>\z@\unskip\fi~\frqq}}

\begin{document}
\newcommand{\Alb}{{\rm Alb}}

\newcommand{\Jac}{{\rm Jac}}

\newcommand{\Hom}{{\rm Hom}}

\newcommand{\End}{{\rm End}}

\newcommand{\Aut}{{\rm Aut}}

\newcommand{\NS}{{\rm NS}}

\title{Elliptic curve configurations on Fano surfaces}

\author{Xavier Roulleau}
\begin{abstract}
The elliptic curves on a surface of general type constitute an obstruction
for the cotangent sheaf to be ample. In this paper, we give the classification
of the configurations of the elliptic curves on the Fano surface of
a smooth cubic threefold. That means that we give the number of such
curves, their intersections and a plane model. This classification
is linked to the classification of the automorphism groups of theses
surfaces. 
\end{abstract}
\maketitle


\section{Introduction.}

Considering a variety $S$, it is a natural question to study the
ampleness of its cotangent sheaf. Recall that, by definition, a bundle
$\mathcal{E}$ on $S$ is ample if the tautological sheaf $\mathcal{O}_{\mathbb{P}(\mathcal{E}^{*})}(1)$
is an ample line bundle of the projective bundle $\mathbb{P}(\mathcal{E}^{*})$
of one dimensional sub-spaces of $\mathcal{E}^{*}$. Gieseker gives
the following criteria of ampleness if in addition the bundle $\mathcal{E}$
is generated by its space of global sections:
\begin{prop}
(Gieseker). The bundle $\mathcal{E}$ is ample if and only if for
every curve $C\hookrightarrow S$, the bundle $\mathcal{E}\otimes\mathcal{O}_{C}$
has no quotient isomorphic to the trivial sheaf $\mathcal{O}_{C}$.
\end{prop}
Hence, we consider varieties with the following assumptions:
\begin{condition}
\label{l'Hypoth=0000E8se}The variety $S$ is a smooth complex surface
of general type. The cotangent sheaf $\Omega_{S}$ of $S$ is generated
by its global sections and the irregularity $q=\dim H^{0}(\Omega_{S})$
satisfies $q>3$.
\end{condition}
With the criteria of Gieseker in mind, a curve $C\hookrightarrow S$
is called non-ample if and only if the sheaf $\Omega_{S}\otimes\mathcal{O}_{C}$
has a quotient isomorphic to the trivial sheaf $\mathcal{O}_{C}$.
These curves are the obstruction to the ampleness of $\Omega_{S}$.
For example, a smooth curve of genus $1$ on $S$ is a non-ample curve.\\
Let $T_{S}=\Omega_{S}^{*}$ be the tangent sheaf and let $\pi:\mathbb{P}(T_{S})\rightarrow S$
be the projection. As $\pi_{*}\mathcal{O}_{\mathbb{P}(T_{S})}(1)=\Omega_{S}$
and the cotangent sheaf is generated by its global sections, we can
define a map: \[
\psi:\mathbb{P}(T_{S})\rightarrow\mathbb{P}(H^{0}(\Omega_{S})^{*})=\mathbb{P}^{q-1}\]
called the cotangent map of $S$, such that $\mathcal{O}_{\mathbb{P}(T_{S})}(1)=\psi^{*}(\mathcal{O}_{\mathbb{P}^{q-1}}(1)).$

This map is the object of the paper \cite{Roulleau}. In the present
paper, we study the cotangent map and the configurations of the non-ample
curves on Fano surfaces. 

These surfaces were discovered by G. Fano and interest in them has
been stimulated by the works of H. Clemens, P. Griffiths \cite{Clemens}
and A. Tyurin \cite{Tyurin}, \cite{Tyurin1} in $1971$. \\
By definition, a Fano surface is the Hilbert scheme of lines of
a smooth cubic threefold $F\hookrightarrow\mathbb{P}^{4}$. This scheme
is a surface $S$ that verifies Hypothesis \ref{l'Hypoth=0000E8se}
and has irregularity $q=5$.\\
By \cite[Tangent Bundle Theorem 12.37]{Clemens}, the image of
the cotangent map $\psi:\mathbb{P}(T_{S})\rightarrow\mathbb{P}(H^{0}(\Omega_{S})^{*})$
of $S$ is a hypersurface $F'$ of $\mathbb{P}(H^{0}(\Omega_{S})^{*})\simeq\mathbb{P}^{4}$
that is isomorphic to the original cubic $F$. Moreover, when we identify
$F$ and $F'$, the triple $(\mathbb{P}(T_{S}),\pi,\psi)$ is the
universal family of lines on $F$.

The main results presented here are:
\begin{thm}
\label{thm:The-non-ample}A) The non-ample curves of a Fano surface
$S$ are the smooth curves of genus $1$.\\
B) There are only $10$ configurations of smooth curves of genus
$1$ on the Fano surfaces. We know a plane model of each curve, the
number of these curves and their intersection numbers.\\
C) For each Fano surface $S$, we construct a particular sub-group
${\bf G}_{S}$ of its automorphism group. This construction goes in
such a way that the knowledge of the group ${\bf G}_{S}$ gives the
knowledge of the elliptic curve configurations on $S$ and reciprocally,
the knowledge of the elliptic curve configuration on $S$ determines
the group ${\bf G}_{S}$.
\end{thm}
The results are summarized in the Classification Theorem \ref{THE CLASSIFICATION THEOREM}.
The number $n_{S}$ of elliptic curves on a Fano surface $S$ verifies
$0\leq n_{S}\leq30$. The surfaces for which $n_{S}>0$ constitute
a $7$ dimensional family in the $10$ dimensional moduli space of
Fano surfaces. \\
The Picard number $\rho_{S}$ of a Fano surface $S$ satisfies
$1\leq\rho_{S}\leq25$ and we can prove that it is $1$ for $S$ generic.
In fact, the number of elliptic curves is linked with the Picard number
$\rho_{S}$: we have $\rho_{S}\geq n_{S}$ unless $n_{S}=30$, in
which case $\rho_{S}=25$. The number $n_{S}$ is also the number
of $1$ dimensional fibres of the cotangent map $\psi$. This shows
that the geometric properties of $\psi$ and the ampleness of the
cotangent bundle vary non-trivially with the Fano surface.\\
We remark also that the results A) and B) of Theorem \ref{thm:The-non-ample}
are analogous to the classical statement on the canonical bundle of
a minimal surface of general type which is ample if and only if the
surface does not contain a $(-2)$-curve and the classification of
canonical surfaces singularities.

We end this paper by an application of our study of Fano surfaces
to construct cubic threefold whose intermediate Jacobian is isomorphic
to a product of elliptic curves as an Abelian variety. The interest
of this results is that in order to prove that the cubic threefolds
are not rational, Clemens and Griffiths use the fact that their intermediate
Jacobian cannot be isomorphic to a product of Jacobians of curves
as a principaly polarized Abelian variety.

\thanks{Part of this paper was writted at the Max-Plank Institute of Bonn,
which is gratefully acknowledged.}

\section{Properties of the Fano surfaces and of their cotangent map.}

\subsection{Properties of the cotangent map.\label{parag 1.1}}

Let $S$ be a surface which verifies Hypothesis \ref{l'Hypoth=0000E8se}.
In the introduction, we defined the cotangent map \[
\psi:\mathbb{P}(T_{S})\rightarrow\mathbb{P}(H^{0}(\Omega_{S})^{*})=\mathbb{P}^{q-1}\]
by the surjective morphism $H^{0}(\Omega_{S})\otimes O_{\mathbb{P}(T_{S})}\rightarrow O_{\mathbb{P}(T_{S})}(1)$.
Here, we state general results about this morphism which will be used
in the sequel ; a complete treatment can be found in \cite{Roulleau}. 

Let us recall that $\pi:\mathbb{P}(T_{S})\rightarrow S$ is the projection.
Let $s$ be a point of $S$. The restriction of the invertible sheaf
$\mathcal{O}_{\mathbb{P}(T_{S})}(1)$ to the fibre $\pi^{-1}(s)\simeq\mathbb{P}^{1}$
is the degree $1$ invertible sheaf and the image under $\psi$ of
that fibre is a line ; we denote this line by $L_{s}\hookrightarrow\mathbb{P}^{q-1}$.

Let $G(2,q)=G(2,H^{0}(\Omega_{S})^{*})$ be the Grassmannian of projective
lines in $\mathbb{P}^{q-1}$. By the universal property of the Grassmannian,
the cotangent map induces a map \[
\mathcal{G}:S\rightarrow G(2,q)\]
called the Gauss map of $S$. The image by $\mathcal{G}$ of a point
$s$ in $S$ is the line $L_{s}$.

Let $A$ be the Albanese variety $S$ ; its tangent space at $0$
is $H^{0}(\Omega_{S})^{*}$. 
\begin{prop}
\label{l'espace sous jacent tangent ..c elll..} Let $C\hookrightarrow S$
be a curve. The image under $\psi$ of $\pi^{-1}(C)$ is a cone if
and only if $C$ is a non-ample curve. \\
An elliptic curve $E\hookrightarrow S$ is a non-ample curve. Let
$p_{E}$ be the vertex of the cone $\psi(\pi^{-1}(E))$. The underlying
space of the point $p_{E}$ is the tangent space to the elliptic curve
$\vartheta(E)$ translated in $0$.\end{prop}
\begin{proof}
The first assertion is \cite[lemme 2.1]{Roulleau}. \\
Let $E$ be an elliptic curve on $S$. The natural quotient $\Omega_{S}\otimes\mathcal{O}_{E}\rightarrow\Omega_{E}$
is a trivial quotient, hence $E$ is non-ample. The last assertion
results from the definition of the cotangent map.\end{proof}
\begin{prop}
\label{Proposition image smooth implies elliptic curve or fiber}\cite[cor. 2.12]{Roulleau}
Let $C\hookrightarrow S$ be a non-ample curve on a surface with irregularity
$q>4$. Suppose that the vertex $p$ of the cone $T=\psi(\pi^{-1}(C))$
is a smooth point of the image of $\psi$. Then $T$ is contained
in the projective tangent space at $p$ of $F\hookrightarrow\mathbb{P}^{q-1}$
and one of the two following possibilities occurs:\\
a) $C^{2}<0$ and $C$ is a smooth curve of genus $1$.\\
b) $C^{2}=0$ and an integral multiple of $C$ is a fiber of a
fibration $f:S\rightarrow B$ onto a curve of genus $b$ with $q-3\leq b\leq q-2$.
\end{prop}
Let $C\hookrightarrow S$ be a curve and let $K$ be a canonical divisor
of S.
\begin{prop}
\label{intersection KC}\cite[prop. 1.20]{Roulleau} The degree of
the cycle $\psi_{*}\pi^{*}C$ equals $KC$.
\end{prop}

\subsection{Main properties of Fano surfaces.}

We recall here some properties of Fano surfaces needed in the sequel.
The main references are the works of Clemens, Griffiths \cite{Clemens}
and of Tyurin \cite{Tyurin}.

Let $F$ be a smooth cubic hypersurface of $\mathbb{P}^{4}$ and let
$S$ be its Fano surface of lines. This surface verifies Hypothesis
\ref{l'Hypoth=0000E8se} and has irregularity $q=5$. Moreover:
\begin{thm}
\label{Tangent Bundle Theorem}\cite[Tangent Bundle Theorem 12.37]{Clemens}
The image of the cotangent map of $S$ is a cubic hypersurface $F'\hookrightarrow\mathbb{P}(H^{0}(\Omega_{S})^{*})\simeq\mathbb{P}^{4}$
isomorphic to the cubic $F\hookrightarrow\mathbb{P}^{4}$. Under the
identification of $H^{0}(\Omega_{S})^{*}$ and $H^{0}(\mathbb{P}^{4},\mathcal{O}_{\mathbb{P}^{4}}(1))$,
and $F'$ and $F$, the triple $(\mathbb{P}(T_{S}),\pi,\psi)$ is
the universal family of lines of $F$. 
\end{thm}
The Chern numbers of a Fano surface verify : $c_{1}^{2}=45,\, c_{2}=27$.
The cotangent map has degree $6$: there are $6$ lines through a
generic point of $F$.
\begin{lem}
\cite[Thm. 12.37]{Clemens}, \cite[Cor. Parag. 4]{Beauville} The
Gauss map and the Albanese map are embeddings. 
\end{lem}
Let $s$ be a point of $S$ and let $C_{s}$ be the closure of points
$t\not=s$ in $S$ such that the line $L_{t}$ intersects the line
$L_{s}$.
\begin{prop}
\label{trois fois le diviseur incident est}\cite[Parag. 10]{Clemens}
The incidence divisor $C_{s}$ is ample, has self-intersection $C_{s}^{2}=5$
and arithmetical genus $11$. The divisor $3C_{s}$ is numerically
equivalent to a canonical divisor.
\end{prop}

\subsection{Properties of a non-ample curve on a Fano surface.}

Let $S$ be a Fano surface and let $F\hookrightarrow\mathbb{P}^{4}$
be the image of its cotangent map $\psi$. For a point $p$ of $F$,
we denote by $T_{F,p}\hookrightarrow\mathbb{P}^{4}$ the projective
tangent hyperplane to $F$ at $p$. For $E\hookrightarrow S$ a non-ample
curve, we denote by $p_{E}$ the vertex of the cone $\psi(\pi^{-1}(E))$.
\begin{prop}
\label{bijection courbe ell et cones}\label{que des courbes elliptiques}
A curve $E$ on $S$ is non-ample if and only if it is an elliptic
curve. \\
Let $E\hookrightarrow S$ be elliptic curve. The cone $\psi(\pi^{-1}(E))$
is the section of $F$ by the hyperplane $T_{F,p_{E}}$, furthermore
: $E^{2}=-3,\textrm{ }C_{s}E=1$.\end{prop}
\begin{proof}
Let $E\hookrightarrow S$ be a non-ample curve. By Proposition \ref{Proposition image smooth implies elliptic curve or fiber},
the intersection of $F$ by $T_{F,p_{E}}$ contains the cone $\psi(\pi^{-1}(E))$.
As a smooth cubic threefold does not contain a plane, the hyperplane
section of $F$ by $T_{F,p_{E}}$ is irreducible and is the cone $\psi(\pi^{-1}(E))$
. This cone has degree $3$. As the Gauss map is an embedding, the
restriction of $\psi$ on $\pi^{-1}(E)$ is birational onto its image,
thus $\psi_{*}\pi^{*}E=\psi(\pi^{-1}(E))$ and Proposition \ref{intersection KC}
implies that \[
KE=\deg\psi_{*}\pi^{*}E=3,\]
where $K$ is a canonical divisor of $S$.\\
By Proposition \ref{Proposition image smooth implies elliptic curve or fiber},
the curve $E$ satifies $E^{2}\leq0$. The number $2p_{a}(E)-2=E^{2}+KE=E^{2}+3$
must be divisible by $2$, hence $E^{2}\not=0$. Proposition \ref{Proposition image smooth implies elliptic curve or fiber}
imply that $E$ is an elliptic curve, thus: $E^{2}=-3$. Since $K$
is numerically equivalent to $3C_{s}$, we obtain: $EC_{s}=1$.
\end{proof}
The following proposition is \cite[Parag. 8, 10]{Clemens}:
\begin{prop}
\label{30courbes} A Fano surface contains at most $30$ smooth curves
of genus $1$. 
\end{prop}

\subsection{The automorphism groups of the cubic and of the Fano surface.\label{parag automorphismes}}

Let us denote by $\Aut(X)$ the automorphisms group of a variety $X$.
Let $S$ be a Fano surface and let $F$ be the image of the cotangent
map of $S$. \\
An element $h\in\Aut(F)$ preserves the lines and by the Tangent
Bundle Theorem \ref{Tangent Bundle Theorem}, it acts on $S$ by an
element denoted by $\rho(h)$.\\
Let be $\tau\in\Aut(S)$. The automorphism $\tau$ acts on the
space $H^{0}(\Omega_{S})^{*}$ by an automorphism denoted by $\tau^{*}\in GL(H^{0}(\Omega_{S})^{*})$.
We remark also that $\tau$ induce an automorphism of the Albanese
variety of $S$ : the morphism $\tau^{*}$ is the differential of
that automorphism.\\
We denote by $\widetilde{\tau}\in PGL(H^{0}(\Omega_{S})^{*})$
the projectivisation of $\tau^{*}$. 

The following proposition is an immediate consequence of the definitions
and the Tangent Bundle Theorem \ref{Tangent Bundle Theorem} ; we
skip its proof because of its length.
\begin{prop}
\label{isomorphisme entre aut S et aut F}The morphism $\rho:\Aut(F)\rightarrow\Aut(S)$
is an isomorphism and its inverse is the morphism $\Aut(S)\rightarrow\Aut(F);\,\tau\rightarrow\widetilde{\tau}$.
\end{prop}
In particular, for a point $s$ of $S$ and $\tau\in\Aut(S)$, we
have : $L_{\tau s}=\widetilde{\tau}(L_{s})$.

\section{Configurations of the elliptic curves.}

\subsection{Configurations of $2$ or $3$ elliptic curves.\label{s poss=0000E8de une c elliptique}}

For two sub-varieties $V_{1},V_{2}$ of $\mathbb{P}^{4}$, we denote
by $\left\langle V_{1},V_{2}\right\rangle $ their linear hull. \\
Let $E\hookrightarrow S$ be an elliptic curve on the Fano surface
$S$ and let $p_{E}$ be the vertex of the cone $\psi(\pi^{-1}(E))$
: this cone is the intersection of $F$ and $T_{F,p_{E}}$. \\
Let $s$ be a point of $S$ outside the curve $E$ : the line $L_{s}$
corresponding to the point $s$ is not inside the cone $\psi(\pi^{-1}(E))$
and $X_{s}:=\left\langle L_{s},p_{E}\right\rangle $ is a plane. This
plane cuts the cubic $F$ in three lines:\\
1) the line $L_{s}$, \\
2) the line $L_{\gamma_{E}s}$ (on the cone) through the vertex
$p_{E}$ and the intersection point of $L_{s}$ and the hyperplane
$T_{F,p_{E}}$, \\
3) the residual line $L_{\sigma_{E}s}$ such that : \[
X_{s}F=L_{s}+L_{\gamma_{E}s}+L_{\sigma_{E}s}.\]
As an Albanese morphism of $S$ is an embedding, the surface $S$
does not contain a rational curve ; furthermore $E$ has genus $>0$.
Then \cite[Cor. 1.44]{Debarre} implies that the rational maps $\sigma_{E}:S\rightarrow S$
and $\gamma_{E}:S\rightarrow E$ are everywhere defined. As the plane
$X_{s}$ ($s$ in $S\setminus E$) is equal to $X_{\sigma_{E}s}$,
the morphism $\sigma_{E}^{2}$ is the identity on $S\setminus E$,
thus $\sigma_{E}$ is an involutive automorphism. \\
Let $s,t$ be two points of the curve $E\hookrightarrow S$. The
line $L_{s}$ cuts the line $L_{t}$ at the vertex of the cone $\psi(\pi^{-1}(E))$:
thus the point $s$ lies on the incidence divisor $C_{t}$ and there
exists a residual divisor $R_{t}$ such that: \[
C_{t}=E+R_{t}.\]

\begin{thm}
\label{la fibration}a) Let $t$ be a point of $E\hookrightarrow S$.
The divisor $R_{t}$ is the fibre at $t$ of $\gamma_{E}$ and has
arithmetical genus $7$. The morphisms $\gamma_{E}$ and $\sigma_{E}$
 satisfy $\gamma_{E}\sigma_{E}=\gamma_{E}$.\\
b) The automorphism $-\sigma_{E}^{*}\in GL(H^{o}(\Omega_{S})^{*})$
is a complex order $2$ reflection. The eigenspace of $\sigma_{E}^{*}$
with eigenvalue $1$ is the underlying space of the vertex $p_{E}\in\mathbb{P}^{4}$.
\\
c) Let $E'\hookrightarrow S$ be an elliptic curve, distinct from
$E$. Then $0\leq EE'\leq1$.\\
d) The automorphisms $\sigma_{E}$ and $\sigma_{E'}$ verify: \[
(\sigma_{E}\sigma_{E'})^{3-EE'}=Id_{S}.\]
e) If $EE'=1$, then the fibration $\gamma_{E}$ contracts $E'$.
\\
f) If $EE'=0$, then there exists a third elliptic curve $E''$
such that \[
\sigma_{E}(E')=\sigma_{E'}(E)=E''.\]
 Moreover the curves $E'$ and $E''$ are sections of $\gamma_{E}$. 
\end{thm}
Let us prove Theorem \ref{la fibration}. Let $t$ be a point of $E\hookrightarrow S$.
Let $s$ be a generic point of $R_{t}=C_{t}-E$. The line $L_{s}$
cuts the line $L_{t}\hookrightarrow\psi(\pi^{-1}(E))$ in a point
different from $p_{E}$ and by definition: $\gamma_{E}s=t$, thus
$s$ is a point of $\gamma_{E}^{-1}(t)$. This proves that $R_{t}$
is a component of $\gamma_{E}^{-1}(t)$. Conversely, we have $R_{t}^{2}=(C_{t}-E)^{2}=5-2+(-3)=0$,
by \cite[Chap. III, Zariski's Lemma 8.2]{Barth}, that implies that
$R_{t}$ is the fibre at $t$ of $\gamma_{E}$.\\
A canonical divisor $K$ is numerically equivalent to $3C_{t}$,
hence $KR_{t}=K(C_{t}-E)=12$ and $R_{t}$ has arithmetical genus
$\frac{12+0}{2}+1$.\\
 The plane $X_{s}$ ($s$ in $S\setminus E$) is equal to the plane
$X_{\sigma_{E}s}$, hence $\gamma_{E}s=\gamma_{E}\sigma_{E}s$.\\
Let $E'\hookrightarrow S$ be an elliptic curve. As $R_{t}$ is
a fibre, we have: \[
R_{t}E'=(C_{t}-E)E'=1-EE'\geq0.\]
If $E=E'$, we see that $\gamma_{E}$ has degree $4$ on $E$. Suppose
now that $E\not=E'$, then $EE'\geq0$ and hence $0\leq EE'\leq1$.
We proved a) and c).

If $EE'=1$, then $R_{t}E=0$ and $E'$ is contained in a fibre of
$\gamma_{E}$, hence e). If $EE'=0$, then $R_{t}E'=1$, thus $E'$
is a section of $\gamma_{E}$. 

Let us prove the property b). Up to a variable change, we can suppose
that the vertex $p_{E}$ of the cone $\psi(\pi^{-1}(E))$ is the point
$(1:0:0:0:0)$. It is easy to prove that, up to a variable change,
an equation of $F$ is :\[
F=\{x_{1}^{2}x_{2}+G(x_{2},x_{3},x_{4},x_{5})=0\}.\]
The automorphism $h:x\rightarrow({-x}_{1}:x_{2}:x_{3}:x_{4}:x_{5})$
acts on $F$. The geometric interpretation of $h$ is given as follows:
if $q$ is a generic point of $F$, then the points $p,q,h(q)$ lie
on a line. Let $s$ be a generic point of $S$. It is then easy to
check that the lines $L_{s}$, $h(L_{s})$ and the point $p_{E}$
span a plane : this is thus the plane $X_{s}$ and we see that $h(L_{s})=L_{\sigma_{E}s}$,
thus : $h=\widetilde{\sigma_{E}}$ (see the notations of Theorem \ref{isomorphisme entre aut S et aut F}).\\
The automorphism $\widetilde{\sigma_{E}}=h$ fix the point $p_{E}$
and the hyperplane $\{x_{1}=0\}\hookrightarrow\mathbb{P}^{4}$ . The
fixed locus of $\sigma_{E}$ is thus the union of $E$ and the $27$
points corresponding to the lines in the intersection of $F$ and
the hyperplane $\{x_{1}=0\}$ (as we can verify this intersection
is a smooth cubic surface). \\
By Proposition \ref{isomorphisme entre aut S et aut F}, the automorphism
$\sigma_{E}^{*}\in GL(H^{0}(\Omega_{S})^{*})$ is equal to: \[
g:x\rightarrow({-x}_{1},x_{2},x_{3},x_{4},x_{5})\]
or to $-g$. Let $\vartheta:S\rightarrow A$ be an Albanese morphism.
As $\sigma_{E}$ is the identity on $E$, the automorphism $\sigma_{E}^{*}\in GL(H^{0}(\Omega_{S})^{*})$
(which is the differential of the action of $\sigma_{E}$ on $A$)
is the identity on the tangent space of the curve $\vartheta(E)$
translated at $0$. By Proposition \ref{l'espace sous jacent tangent ..c elll..},
this space is $\mathbb{C}(1,0,0,0,0)\subset H^{0}(\Omega_{S})^{*}$.
Thus $\sigma_{E}^{*}=-g$ and $-\sigma_{E}^{*}$ is a complex reflection
of order $2$. We proved b).

Let $E'\not=E$ a second smooth curve of genus $1$ on $S$. Let us
prove that $(\sigma_{E}\sigma_{E'})^{3-EE'}=Id_{S}$.

\textbf{Case $EE'=1$.} Suppose that $EE'=1$. A generic hyperplane
section $Y$ of the cone $\psi(\pi^{-1}(E'))$ parametrizes the lines
of this cone and is a plane model of $E'$ : we will identify $Y$
and $E'$. Let $t$ be the intersection point of $E$ and $E'$ :
with this neutral element, the curve $E'$ is an elliptic curve.\\
 Let $s$ be a point of $E'$ different from $t$, then by definition
of $\sigma_{E}$ and $\gamma_{E}$: \[
X_{E,s}F=L_{s}+L_{\gamma_{E}s}+L_{\sigma_{E}s}\]
where $X_{E,s}:=\left\langle L_{s},p_{E}\right\rangle $. Since the
line $L_{s}$ cuts the line $L_{t}$ at the vertex $p_{E'}$, we have
$X_{E,s}=\left\langle L_{s},L_{t}\right\rangle $. Thus the point
$t$ is one of the three points $s,\,\gamma_{E}s,\,\sigma_{E}s$.
As $s\not=t$ and the automorphism $\sigma_{E}$ preserves $E$, the
point $\sigma_{E}s$ is not an element of $E$, hence $\gamma_{E}s=t$.
Thus $E'$ is a component of $\gamma_{E}^{-1}(t)$. \\
The plane $X_{E,s}$ contains the lines $L_{s}$, $L_{t}$ that
are in the hyperplane section $\psi(\pi^{-1}(E'))$ of $F$, thus
the third line $L_{\sigma_{E}s}$ is also in the cone $\psi(\pi^{-1}(E'))$
and the point $\sigma_{E}s$ is on $E'$. Because of the relation
\[
X_{E,s}F=L_{s}+L_{\gamma_{E}s}+L_{\sigma_{E}s},\]
the three points $s,t=\gamma_{E}s,\sigma_{E}s$ are on a line in the
plane model $Y$, hence $\sigma_{E}s=-s$ for all points $s$ of $E'$.
\begin{rem}
\label{Remarque le point p_E' est dans l'hyper. des points fixes}
Since $\sigma_{E}(s)=-s$ on $E'$, the endomorphism $\sigma_{E}^{*}\in GL(H^{0}(\Omega_{S})^{*})$
is the morphism of multiplication by $-1$ on the tangent space to
the curve $\vartheta(E')$ (translated to $0$). Thus : $\widetilde{\sigma_{E}}p_{E'}=p_{E'}$
and $p_{E'}$ is the intersection point of the line $L_{t}$ and the
hyperplane of fixed points of $\widetilde{\sigma_{E}}$. This implies
that the points $p_{E}$ and $p_{E'}$ are the only vertices of cones
on the line $L_{t}$.
\end{rem}
Let $s$ be a generic point of $S$, then\[
\widetilde{\sigma_{E'}}X_{E,s}=\widetilde{\sigma_{E'}}\left\langle p_{E},L_{s}\right\rangle =\left\langle p_{E},L_{\sigma_{E'}s}\right\rangle =X_{E,\sigma_{E'}s},\]
hence: \[
\widetilde{\sigma_{E'}}X_{E,s}F=L_{\sigma_{E'}s}+L_{\gamma_{E}\sigma_{E'}s}+L_{\sigma_{E}\sigma_{E'}s}.\]
But $X_{E,s}F=L_{s}+L_{\gamma_{E}s}+L_{\sigma_{E}s}$ hence: \[
\widetilde{\sigma_{E'}}X_{E,s}F=L_{\sigma_{E'}s}+L_{\sigma_{E'}\gamma_{E}s}+L_{\sigma_{E'}\sigma_{E}s}.\]
Since $\gamma_{E}\sigma_{E'}s$ and $\sigma_{E'}\gamma_{E}s$ are
points of $E$, we see that $\sigma_{E}\sigma_{E'}s=\sigma_{E'}\sigma_{E}s$,
thus $(\sigma_{E}\sigma_{E'})^{2}=Id_{S}$. 

\textbf{Case $EE'=0$.} Suppose now that $EE'=0$. Let $s$ be a point
of $E'$. By definition $X_{E,s}=\left\langle p_{E},L_{s}\right\rangle $
and: \[
X_{E,s}F=L_{s}+L_{\gamma_{E}s}+L_{\sigma_{E}s}.\]
Suppose that $\sigma_{E}s$ is a point of $E'$, then plane $X_{E,s}$
cuts the cone $\psi(\pi^{-1}(E'))$ into two lines : $L_{s}$ and
$L_{\sigma_{E}s}$. As this cone has degree $3$, the intersection
of $X_{E,s}$ and $\psi(\pi^{-1}(E'))$ contains the third line $L_{\gamma_{E}s}$.
Since $\psi(\pi^{-1}(E'))$ is a cone, that implies that $\gamma_{E}s$
is a point of $E'$ : this is a impossible because $EE'=0$. Thus
$\sigma_{E}s$ is not a point of $E'$.\\
 The automorphism $\sigma_{E}$ fixes $E$. Hence the point $\sigma_{E}s$
is not a point of $E$. This proves that the surface $S$ contains
a third smooth curve $E''=\sigma_{E}(E')$ of genus $1$ and that
$EE''=E'E''=0$.\\
For a point $s$ of $E'$, the plane $X_{E,s}=\left\langle p_{E},L_{s}\right\rangle $
contains the line $L_{\gamma_{E}s}$ and the point $p_{E'}\in L_{s}$,
hence: \[
X_{E,s}=\left\langle p_{E},L_{s}\right\rangle =\left\langle p_{E'},L_{\gamma_{E}s}\right\rangle =X_{E',\gamma_{E}s}.\]
But we have \[
X_{E,s}F=L_{s}+L_{\gamma_{E}s}+L_{\sigma_{E}s}\]
 and \[
X_{E',\gamma_{E}s}F=L_{\gamma_{E}s}+L_{\gamma_{E'}\gamma_{E}s}+L_{\sigma_{E'}\gamma_{E}s}.\]
Since the points $s,\,\gamma_{E}s$ and $\sigma_{E}s$ are respectively
points of $E',E$ and $E''$, we see that for all points $s$ of $E'$:
\[
\sigma_{E'}\gamma_{E}s=\sigma_{E}s.\]
Hence the restriction of $\sigma_{E'}$ to $E$ is a morphism from
$E$ to $E''$ and $\sigma_{E'}(E)=\sigma_{E}(E')=E''$, this proves
f).\\
Let $s$ be a point of $E'$. The lines $L_{s},L_{\gamma_{E}s}$
and $L_{\sigma_{E}s}$ contain respectively the vertices $p_{E'},p_{E}$,
and $p_{E''}$. As $s$ varies in $E'$, the plane $X_{E,s}$ varies
and the linear hull $\ell$ of the points $p_{E},\, p_{E'}$ and $p_{E''}$
cannot be a plane: it is a line. 
\begin{rem}
\label{remarque au plus trois points}The line $\ell$ lies outside
the cubic $F$, otherwise the curves $E$ and $E'$ would have a common
point. The points $p_{E},\, p_{E'}$ and $p_{E''}$ are the intersection
points of $\ell$ and the cubic $F$.
\end{rem}
Let $s$ be a point of $S$. The morphism $\widetilde{\sigma}_{E'}$
verifies : $\widetilde{\sigma}_{E'}(L_{s})=L_{\sigma_{E'}s}$, and
furthermore: \begin{equation}
\widetilde{\sigma}_{E'}X_{E,\sigma_{E'}s}=\widetilde{\sigma}_{E'}\left\langle p_{E},L_{\sigma_{E'}s}\right\rangle =\left\langle p_{E''},L_{s}\right\rangle =X_{E'',s}.\label{etiquette1}\end{equation}
We have $X_{E,\sigma_{E'}s}F=L_{\sigma_{E'}s}+L_{\gamma_{E}\sigma_{E'}s}+L_{\sigma_{E}\sigma_{E'}s}$.
Hence \[
\widetilde{\sigma}_{E'}X_{E,\sigma_{E'}s}F=L_{\sigma_{E'}^{2}s}+L_{\sigma_{E'}\gamma_{E}\sigma_{E'}s}+L_{\sigma_{E'}\sigma_{E}\sigma_{E'}s},\]
 but by \ref{etiquette1}: \[
\widetilde{\sigma}_{E'}X_{E,\sigma_{E'}s}F=X_{E'',s}F.\]
Since $X_{E'',s}F=L_{s}+L_{\gamma_{E''}s}+L_{\sigma_{E''}s}$, we
see that $\sigma_{E'}\sigma_{E}\sigma_{E'}=\sigma_{E''}$. So the
group generated by $\sigma_{E},\sigma_{E'},\sigma_{E''}$ is isomorphic
to $\Sigma_{3}$ and $(\sigma_{E}\sigma_{E'})^{3}=Id_{S}$. 

This ends the proof of Theorem \ref{la fibration}. $\Box$

Let us now study the configuration of $3$ elliptic curves.
\begin{prop}
\label{premier lemme}Let $E_{1},E_{2}$ and $E$ be three elliptic
curves on $S$ such that $E_{1}E=E_{2}E=1$. Then: $E_{1}E_{2}=0$
and the curve $E_{3}=\sigma_{E_{1}}(E_{2})$ verifies $E_{3}E=1$.\end{prop}
\begin{proof}
Suppose $E_{1}E_{2}=1$. This implies that the line through $p_{E_{1}}$
and $p_{E_{2}}$ lies on the cone $\psi(\pi^{-1}(E))$ and hence goes
through $p_{E}$. But by remark \ref{Remarque le point p_E' est dans l'hyper. des points fixes},
the three points $p_{E},\, p_{E_{1}}$ and $p_{E_{2}}$ cannot be
on a line. Hence $E_{1}E_{2}=0$.\\
Since $EE_{1}=1$, we have $\sigma_{E_{1}}(E)=E$ and $E_{3}E=\sigma_{E_{1}}(E_{2})\sigma_{E_{1}}(E)=E_{2}E=1$.
\end{proof}

\subsection{The graph of the configuration of the vertices of cones.}

Let $\mathcal{E}$ be the set of elliptic curves on $S$. Let us consider
the following graph $\mathbb{G}$: the set of vertices of $\mathbb{G}$
is $\mathcal{E}$ and an edge links $E\in\mathcal{E}$ to $E'\in\mathcal{E}$
if and only if $EE'=0$. \\
Let ${\bf G}_{S}$ be the sub-group of $\Aut(S)$ generated by
the automorphisms $\sigma_{E},\, E\in\mathcal{E}$. We have the three
following relations between its generators:\\
a) for all $E\in\mathcal{E}$, $\sigma_{E}$ verifies $\sigma_{E}^{2}=Id_{S}$,\\
b) an edge links $E$ and $E'$ if and only if $(\sigma_{E}\sigma_{E'})^{3}=Id_{S}$,
\\
c) otherwise $(\sigma_{E}\sigma_{E'})^{2}=Id_{S}$.\\
The following corollary is a consequence of Proposition \ref{premier lemme}. 
\begin{cor}
\label{si trois sommets alors une ar=0000EAte}Let $E_{1},E_{2}$
and $E_{3}$ be three elements of $\mathcal{E}$. At least one edge
links two of the three vertices $E_{1},E_{2},E_{3}$ of the graph
$\mathbb{G}$.\\
There is no sub-group of ${\bf G}_{S}$ generated by some elements
$\sigma_{E},E\in\mathcal{E}$ and isomorphic to $(\mathbb{Z}/2\mathbb{Z})^{3}$.\end{cor}
\begin{proof}
If there are no edges between the vertices $E_{1}$ and $E_{3}$ and
between the vertices $E_{2}$ and $E_{3}$, then $E_{1}E_{3}=E_{2}E_{3}=1$
and the Proposition \ref{premier lemme} implies that $E_{1}E_{2}=0$.
Thus an edge links the vertices $E_{1}$ and $E_{2}$. The second
assertion is a reformulation of the first.
\end{proof}
We remark that if the graph $\mathbb{G}$ has $m$ connected components,
then the group ${\bf G}_{S}$ is the direct product of $m$ sub-groups.
The Corollary \ref{si trois sommets alors une ar=0000EAte} can be
reformulated as follows:
\begin{cor}
\label{graphe connexe}The graph $\mathbb{G}$ is connected or $\mathbb{G}$
has two connected components $\mathbb{G}_{1},\mathbb{G}_{2}$ such
that two different vertices of a component $\mathbb{G}_{i}$ are linked
by an edge.
\end{cor}
Let $\mathcal{E}'$ be a sub-set of $\mathcal{E}$. Let $\mathbb{G}'$
be the graph whose set of vertices is $\mathcal{E}'$ and such that
an edge links two elements of $\mathcal{E}'$ if and only if these
vertices are linked by an edge in $\mathbb{G}$. Suppose that the
three relations a), b) and c) above are the only ones between the
elements $\sigma_{E},\, E\in\mathcal{E}'$. The group generated by
the automorphisms $\sigma_{E},\, E\in\mathcal{E}'$ is then a Weyl
group, which we denote by ${\bf G}_{S}'$.
\begin{cor}
\label{groupe WAn}If the graph $\mathbb{G}'$ is connected and has
$n$ vertices, then $1\leq n\leq4$ and the group ${\bf G}_{S}'$
is the permutation group of the set of $n+1$ elements.\end{cor}
\begin{proof}
By corollary \ref{si trois sommets alors une ar=0000EAte} and the
classification of the Weyl groups \cite{Humphreys}, the graph $\mathbb{G}'$
must be one of the graphs $A_{n},\,1\leq n\leq4$. The Weyl group
$W(A_{n})$ associated to $A_{n}$ is the permutation group of $n+1$
elements.
\end{proof}

\subsection{Restrictions on the complex reflection groups.\label{paragraphe groupe d'automor homographies 2}}

Let us denote by $G_{S}$ the complex reflection group generated by
the order $2$ reflections $-\sigma_{E}^{*}\in GL(H^{0}(\Omega_{S})^{*})$,
$E\hookrightarrow S$ elliptic curve. For the basic properties of
reflection groups see \cite{todd}, \cite{cohen} or \cite{Dolgachev},
here we recall some them:
\begin{defn}
\label{d=0000E9finition de spe de relexion}A reflection in a space
$V$ is a linear transformation of $V$ of finite order with fixed
point set an hyperplane of $V$. A reflection group is a finite group
generated by reflections.\\
Let $G_{1}$ and $G_{2}$ be two reflection groups acting on spaces
$V_{1}$ and $V_{2}$, we say (improperly) that $G_{1}$ is a reflection
sub-group of $G_{2}$, if there exists an injective morphism $G_{1}\rightarrow G_{2}$
of complex reflection groups. In this case, we denote the elements
of $G_{1}$ and $G_{2}$ by the same letters.
\end{defn}
The list of the $37$ types of irreducible reflection groups was compiled
by G. Shephard and J. Todd \cite{todd}. \\
Let us fix some notation. Let $m,\, n>0$ be integers and let $p$
be an integer dividing $m$. We denote by $\Sigma_{n}\subset GL_{n}(\mathbb{C})$
the group of permutation matrices and by $A(m,p,n)\subset GL_{n}(\mathbb{C})$
the group of diagonal matrices $D$ of order $m$ such that $\det(D)^{\frac{m}{p}}=1$.
The type $2$ reflection groups are the groups $G(m,p,n)$ generated
by $A(m,p,n)$ and $\Sigma_{n}$.\\
The type $3$ groups constitute the family $[\,]^{n}$ for $n\in\mathbb{N}\setminus\{0,1\}$
where the $[\,]^{n}$ is the group of morphisms $\mathbb{C}\rightarrow\mathbb{C};\, x\rightarrow\xi x$
, $\xi^{n}=1$.
\begin{thm}
\label{restriction sur les groupes}An irreducible sub-group of $G_{S}$
generated by reflections of order $2$ is isomorphic to one of the
following groups: \[
\{1\},\,[\,]^{2},\, G(3,3,n),\, G(1,1,n)=\Sigma_{n},\,2\leq n\leq5.\]
\end{thm}
\begin{proof}
By the Remarks \ref{Remarque le point p_E' est dans l'hyper. des points fixes}
and \ref{remarque au plus trois points}, a projective line of $\mathbb{P}^{4}$
contains at most $3$ vertices of a cone, hence a $2$ dimensional
reflection sub-group of $G_{S}$ contains at most $3$ reflections
of order $2$. The groups of type $4$ to $22$ are $2$ dimensional
irreducible groups. They either have no or at least $6$ reflections
of order $2$ \cite[table p.395]{cohen} : none of these groups can
be a reflection sub-group of $G_{S}$.

The groups numbered $25,29,31,32,33,34$ possess either no or at least
$40$ reflections of order $2$ \cite[p.412]{cohen}. Since reflections
of order $2$ of $G_{S}$ are in bijection with elliptic curves on
$S$ and a Fano surface contains at most $30$ elliptic curves (Proposition
\ref{30courbes}), none of these groups is a sub-group of $G_{S}$.

The groups $23,28,30,35,36,37$ are real reflection groups \cite[p.299]{todd}
which have been excluded by the Corollary \ref{groupe WAn}. 

The group $26$ has a reflection sub-group isomorphic to the group
number $4$ already excluded \cite[p.302]{todd}. 

The groups $24$ and $27$ have two reflections $R_{1}$ and $R_{2}$
of order $2$ such that $(R_{1}R_{2})^{4}=1$ \cite[p.299]{todd}
and hence they cannot be sub-groups of $G_{S}$.

It thus remains to study the reflection groups of types $1$, $2$
and $3$.

Let $M$ be an irreducible sub-group of $G_{S}$ generated by order
$2$ reflections.\\
Type $3$. It is immediate that there exists $n>1$ such that $[\,]^{n}\simeq M$
if and only if $n=2$.\\
Type $2$ and $1$. The group $G(m,p,n)$ (where $p$ divides $m\in\mathbb{N}^{*}$
and $n>1$) is an $n$ dimensional irreducible reflection group if
and only if $m>1$ and $(m,p,n)\not=(2,2,2)$. The representation
$G(1,1,n)=\Sigma_{n}$ of the permutation group breaks up into a $1$
dimensional trivial representation and an $n-1$ dimensional irreducible
representation $W(A_{n-1})$ called the standard representation. The
groups: \[
W(A_{n}),\, n\in\mathbb{N}^{*}\]
constitute the number $1$ reflection type in the Shephard-Todd classification. 

Let $m,\, p,\, n$ be integers such that $m>0,$ $p$ divides $m$
and $n>1$. Suppose that $M$ is the group $G(m,p,n)$. Theorem \ref{la fibration}
implies that a group generated by two reflections of order $2$ of
$G_{S}$ is the diedral group of order $4$ or $6$. For $m>1$, the
group $G(m,p,n)$ contains a diedral sub-group of order $2m$ generated
by two reflections of order $2$. Thus the integer $m$ is an element
of $\{1,2,3\}$ ; moreover $n\leq5$.\\
 The group $G(2,2,2)$ is not irreducible and is thus excluded.
The groups $G(2,1,n)$ with $n\geq3$ and the groups $G(2,2,n)$ with
$n\geq4$ contain sub-groups isomorphic to $(\mathbb{Z}/2\mathbb{Z})^{3}$
generated by reflections of order $2$. By Corollary \ref{graphe connexe},
such groups cannot be reflection sub-groups of $G_{S}$.\\
The reflection group $\Sigma_{4}$ is isomorphic to the group $G(2,2,3)$
plus the trivial representation. So $G(2,2,3)$ is implicitly in the
list of Theorem \ref{restriction sur les groupes}.\\
The group $G(3,1,n)$ cannot be isomorphic to $M$ because its
order $2$ reflections generate the strict sub-group $G(3,3,n)\subset G(3,1,n)$.\\
By Corollary \ref{groupe WAn}, if $W(A_{n})$ is a reflection
sub-group of $G_{S}$, then $n\leq4$. The group $\Sigma_{n}$ is
isomorphic to the reflection group $W(A_{n-1})$ plus the trivial
representation.

Hence, we proved that $M$ is isomorphic to one of the groups $\{1\},\,[\,]^{2},\,$$G(3,3,n)$,
$G(1,1,n)=\Sigma_{n}$$,\,2\leq n\leq5$.
\end{proof}

\subsection{Classification of groups and Fano surfaces. \label{paragraphe classif des cubiques}}

Let us classify the Fano surfaces according to the configuration of
their elliptic curves. Let us recall the notations : to each elliptic
curve $E\hookrightarrow S$ corresponds an automorphism $\sigma_{E}$
of $S$ and the group $G_{S}$ is generated by the elements $-\sigma_{E}^{*}\in GL(H^{0}(\Omega_{S})^{*})$,
where the automorphisms $\sigma_{E}^{*}$ are defined in paragraph
\ref{parag automorphismes}. We will use the following remark:
\begin{rem}
Let $F_{eq}$ be an equation of the image $F$ of the cotangent map
of a Fano surface $S$. There exists a morphism: \[
\chi:G_{S}\rightarrow\mathbb{C}^{*}\]
such that $F_{eq}\circ N=\chi(N)F_{eq}$ for all $N\in G_{S}$. 
\end{rem}
Thus we are looking for cubic forms $F_{eq}$, reflection groups $G$
and morphisms $\chi:G\rightarrow\mathbb{C}^{*}$ such that $F_{eq}\circ N=\chi(N)F_{eq}$
for all $N\in G$ and such that $\{F_{eq}=0\}$ is smooth. 

We need some notations and preliminary materials.\\
The order of $G(m,p,n)$ is equal to $\frac{m}{p}m^{n-1}n!$ and
the number of its order $2$ reflections is $m\frac{n(n-1)}{2}$.
The group $G(m,m,n)$ acts on the polynomial space of $\mathbb{C}^{n}$.
The algebra of invariant polynomials is generated by the polynomials:
\[
\sum_{i=1}^{i=n}x_{i}^{mk},\: k\in\{1,...,n-1\}\]
 and by $x_{1}x_{2}....x_{n}$ (see \cite{todd}). \\
Let $S$ be a Fano surface such that the group $G(m,m,n)$ ($m\in\{1,3\},2\leq n\leq5$)
is a reflection sub-group of $G_{S}$. Let be $m\in\{1,3\}$ and $n>1$.
We easily check that the only non-trivial morphism from $G(m,m,n)$
to $\mathbb{C}^{*}$ is the determinant. We call a polynomial $P$
an anti-invariant of $G(m,m,n)$ if $P\circ N=(\det N)P$ for all
$N\in G(m,m,n)$. We verify that:
\begin{lem}
The only reflection groups $G(m,m,n)$ with $m\in\{1,3\}$ and $n\geq2$
that possess an anti-invariant polynomial of degree $\leq3$ are $G(1,1,2)$,
$G(3,3,2)$ and $G(1,1,3)$.\end{lem}
\begin{notation}
\label{notations base duale}For $\lambda^{3}\not=1$, we denote by
$E_{\lambda}$ the smooth plane cubic: $x^{3}+y^{3}+z^{3}-3\lambda xyz=0.$\\
We denote by $A$ the Albanese variety of $S$, by $\vartheta:S\rightarrow A$
an Albanese morphism, by $e_{1},\dots,e_{5}$ the dual basis of the
basis $x_{1},\dots,x_{5}\in H^{o}(\Omega_{S})$. If $v\in H^{0}(\Omega_{S})^{*}$
is a non zero vector, $\mathbb{C}v$ is the vector space generated
by $v$ or the point of $\mathbb{P}^{4}=\mathbb{P}(H^{0}(\Omega_{S})^{*})$
corresponding to this space, we specify as need be. We denote by $\mu_{3}$
the group of third roots of unity. 
\end{notation}
Recall (see Theorem \ref{la fibration} and its proof) that there
is a one to one correspondance between :\\
a) the elliptic curves $E$ on $S$.\\
b) the order $2$ reflections $R=-\sigma_{E}^{*}$ of $G_{S}$.\\
c) the cones $\psi(\pi^{-1}(E))$ on $F$.\\
d) the vertices $p_{E}$ of the cones $\psi(\pi^{-1}(E))\hookrightarrow F$.
\\
When we consider an order $2$ reflection $R$ corresponding to
an elliptic curve $E\hookrightarrow S$, the vertex $p_{E}$ of the
cone $\psi(\pi^{-1}(E))$ is the point of $\mathbb{P}^{4}$ corresponding
to the eigenspace with eigenvalue $-1$ of $R\in GL(H^{0}(\Omega_{S})^{*})=GL_{5}(\mathbb{C})$.
\\
The cone $\psi(\pi^{-1}(E))$ is the intersection of the cubic
$F$ and the tangent space at $p_{E}$. A plane model of the curve
$E$ is the intersection of this cone with the hyperplane of fixed
points of $\widetilde{\sigma_{E}}$ (where $\widetilde{\sigma_{E}}$
is the projectivization of $R$). \\
In order to know the intersection number between two elliptic curves
$E_{1}$ and $E_{2}$ corresponding to the reflections $R_{1}$ and
$R_{2}$, it suffice to compute the order ( $=2$ or $3$) of the
element $R_{1}R_{2}$ and to use the formula : \[
(R_{1}R_{2})^{3-E_{1}E_{2}}=Id.\]
It can also be verified by hand if the line between the vertices $p_{E_{1}}$
and $p_{E_{2}}$ is or not inside the cubic, accordingly: $E_{1}E_{2}=1$
or $0$.\\

In the sequel, we proceed to the classification of the elliptic curve
configurations on Fano surfaces according to the dimension of the
studied irreducible reflection sub-group of $G_{S}$ ; this group
must be in the list of Theorem \ref{restriction sur les groupes}. 

- 1 dimensional sub-group of $G_{S}$.

We proved in paragraph \ref{s poss=0000E8de une c elliptique}, that
when $G_{S}$ contains the group $[\,]^{2}$ generated by a reflection
$R$, the cubic can be writted as follows :\[
F=\{x_{1}^{2}x_{2}+G(x_{2},..,x_{5})=0\}\]
 where $G$ is a cubic form and $R$ acts by $x\rightarrow(-x_{1},x_{2},x_{3},x_{4},x_{5})$.
The tangent space at $p=(1:0:0:0:0)$ is $T_{F,p}=\{x_{2}=0\}$. This
point $p$ is the vertex of a cone in $F$. This cone is the intersection
of $T_{F,p}$ and $F$. Let us denote by $E\hookrightarrow S$ the
elliptic curve on the Fano surface that parametrizes the lines on
the cone. \\
The hyperplane $\{x_{1}=0\}$ is an invariant of the cubic because
it is the fixed locus of $\widetilde{\sigma_{E}}$. A plane model
of $E$ is obtained by the intersection of $F$ and the plane $\{x_{1}=x_{2}=0\}$.

We see in that way, that given a pair $(T,E)$ where $T$ is a smooth
cubic surface (say given by $T=\{G(x_{2},..,x_{5})=0\}$) and $E$
is a smooth hyperplane section of $T$ (say $E=T\cap\{x_{2}=0\}$),
we can associate the pair $(S,E')$ where $S$ is the Fano surface
of $\{x_{1}^{2}x_{2}+G(x_{2},..,x_{5})=0\}$ and $E'$ is an elliptic
curve on it isomorphic to $E$. Reciprocally, given the pair $(S,E')$,
we can recover the cubic $F$ and the pair $(T,E)$, up to isomorphism.
As the moduli space of isomorphism class of such pairs $(T,E)$ is
$7$ dimensional, the moduli of Fano surfaces that contains an elliptic
curve is a $7$ dimensional sub-space of the $10$ dimensional moduli
of Fano surfaces.

- $2$ dimensional reflection sub-groups of $G_{S}$. 

$\blacksquare$ The anti-invariant polynomials of the reflection group
$G(3,3,2)$ yield singular cubics. The invariants of the group $G(3,3,2)$
are generated by the polynomials: \[
x_{1}^{3}+x_{2}^{3},\, x_{1}x_{2},\, x_{3},x_{4},x_{5}.\]
 Up to a change of variables, the cubic $F$ is: \[
F=\{x_{1}^{3}+x_{2}^{3}+3x_{1}x_{2}l(x_{3},x_{4},x_{5})+x_{3}^{3}+x_{4}^{3}+x_{5}^{3}-3\lambda x_{3}x_{4}x_{5}=0\}\]
 where $l$ is a linear form and $\lambda\in\mathbb{C}$. \\
 The $3$ points $\mathbb{C}(e_{1}-\beta e_{2})$ , $\beta\in\mu_{3}$
(see notations \ref{notations base duale}) are vertices of cones.
The Fano surface $S$ contains three disjoint elliptic curves $E_{12}^{\beta},\,\beta\in\mu_{3}$
that have the following model: \[
\{x_{3}^{3}+x_{4}^{3}+x_{5}^{3}-3\lambda x_{3}x_{4}x_{5}+l^{3}=0\}.\]
 Note that we also have studied the reflection group $G(1,1,3)=\Sigma_{3}$
because its representation decomposes into the trivial one plus the
standard representation $W(A_{2})$ which is equal to $G(3,3,2)$. 

- $3$ dimensional reflection sub-groups of $G_{S}$.

$\blacksquare$ The invariants of $G(3,3,3)$ are generated by:\[
x_{1}^{3}+x_{2}^{3}+x_{3}^{3},\, x_{1}x_{2}x_{3},\, x_{4},\, x_{5}.\]
Up to a variables change, there exist coordinates there exists $\lambda\in\mathbb{C}$
such that: \[
F=\{x_{1}^{3}+x_{2}^{3}+x_{3}^{3}-3\lambda x_{1}x_{2}x_{3}+x_{4}^{3}+x_{5}^{3}=0\}.\]
 The group $G(3,3,3)\times G(3,3,2)$ is a reflection sub-group of
$G_{S}$. The $9$ points $\mathbb{C}(e_{i}-\beta e_{j}),\,1\leq i<j\leq3,\,\beta\in\mu_{3}$
are vertices of cones. The corresponding elliptic curves : \[
E_{ij}^{\beta},\,1\leq i<j\leq3,\,\beta\in\mu_{3}\]
on the Fano surface $S$ are isomorphic to the curve $E_{0}$ and
are disjoint.\\
The $3$ points $\mathbb{C}(e_{4}-\beta e_{5}),\,\beta\in\mu_{3}$
are vertices of cones. The corresponding $3$ elliptic curves $E_{45}^{\beta},\,\beta\in\mu_{3}$
are disjoint and isomorphic to the curve $E_{\lambda}$. For all $1\leq i<j\leq3,\,(\beta,\gamma)\in\mu_{3}^{2}$,
we have: \[
E_{ij}^{\beta}E_{45}^{\gamma}=1.\]
The Fano surface contains $12$ elliptic curves.

$\blacksquare$ The invariants of degree less than or equal to $3$
of $\Sigma_{4}$ are generated by: \[
x_{5}^{k},\; x_{1}^{k}+x_{2}^{k}+x_{3}^{k}+x_{4}^{k},\; k\in\{0,1,2,3\}.\]
 Let $F$ be a smooth cubic defined by an element of the $6$ dimensional
space of invariant cubics. For $1\leq i<j\leq4$, the point $\mathbb{C}(e_{i}-e_{j})$
of $\mathbb{P}^{4}$ is the vertex of a cone on $F$ ; let us denote
by: \[
E_{ij}\hookrightarrow S\]
the corresponding elliptic curve on the Fano surface $S$. We have:
\[
E_{ij}E_{st}=\left\{ \begin{array}{cc}
1 & \textrm{if }\{i,j\}\cap\{s,t\}=\emptyset\\
-3 & \textrm{if }E_{ij}=E_{st}\\
0 & \textrm{otherwise.}\end{array}\right.\]
The surface $S$ contains $6$ elliptic curves.

- $4$ dimensional reflection sub-groups of $G_{S}$.

$\blacksquare$ Suppose that $\Sigma_{5}$ is a reflection sub-group
of $G_{S}$. The invariants of degree less than or equal to $3$ of
$\Sigma_{4}$ are generated by: \[
x_{1}^{k}+x_{2}^{k}+x_{3}^{k}+x_{4}^{k}+x_{5}^{k},\; k\in\{0,1,2,3\}.\]
There exist $\lambda,\mu\in\mathbb{C}$ such that the image of the
cotangent map of $S$ is: \[
F=\{\sum_{i=1}^{i=5}x_{i}^{3}+\lambda(\sum_{i=1}^{i=5}x_{i})(\sum_{i=1}^{i=5}x_{i}^{2})+\mu(\sum_{i=1}^{i=5}x_{i})^{3}=0\}.\]
 For $1\leq i<j\leq5$, the point \[
p_{ij}=\mathbb{C}(e_{i}-e_{j})\]
is the vertex of a cone and we denote by $E_{ij}\hookrightarrow S$
the corresponding elliptic curve. We have: \[
E_{ij}E_{st}=\left\{ \begin{array}{cc}
1 & \textrm{if }\{i,j\}\cap\{s,t\}=\emptyset\\
-3 & \textrm{if }E_{ij}=E_{st}\\
0 & \textrm{otherwise.}\end{array}\right.\]
The dual graph of this configuration of $10$ elliptic curves is the
trivalent Petersen graph.

- $5$ dimensional reflection sub-group of $G_{S}$.

$\blacksquare$ The Fermat cubic:\[
F=\{x_{1}^{3}+x_{2}^{3}+x_{3}^{3}+x_{4}^{3}+x_{5}^{3}=0\}\]
is, up to isomorphism, the only cubic stable under $G(3,3,4)$ and
$G(3,3,5)$. The points $\mathbb{C}(e_{i}-\beta e_{j}),\,1\leq i<j\leq5,\,\beta\in\mu_{3}$
are vertices of cones on $F$. Its Fano surface $S$ is the unique
Fano surface that contains $30$ smooth curves of genus $1$. These
curves are numbered:\[
E_{ij}^{\beta},\,1\leq i<j\leq5,\,\beta\in\mu_{3}.\]
Let $E_{ij}^{\gamma}$ and $E_{st}^{\beta}$ be two such curves, then:
\[
E_{ij}^{\beta}E_{st}^{\gamma}=\left\{ \begin{array}{cc}
1 & \textrm{if }\{i,j\}\cap\{s,t\}=\emptyset\\
-3 & \textrm{if }E_{ij}^{\beta}=E_{st}^{\gamma}\\
0 & \textrm{otherwise.}\end{array}\right.\]
A remarkable property of this surface is that its Néron-Severi group
$NS(S)$ has rank $25=\dim H^{1}(S,\Omega_{S})$ and $NS(S)\otimes\mathbb{Q}$
is generated by the $30$ elliptic curves.

Now, we study the case for which the reflection group $G_{S}$ is
not irreducible. Corollary \ref{graphe connexe} proves that :
\begin{lem}
\label{lemme produit de 2 refl}If the reflection group $G_{S}$ is
not irreducible, it is the direct product of two irreducible reflection
groups $\mathbb{W}_{1}$ and $\mathbb{W}_{2}$ such that if $R_{1}$
and $R_{2}$ are two different reflections of order $2$ of the group
$\mathbb{W}_{i}$, we have : \[
(R_{1}R_{2})^{3}=Id.\]

\end{lem}
The groups with this last property and listed in Theorem \ref{restriction sur les groupes}
are: \[
[\,]^{2},\, G(3,3,2)\textrm{ or }G(3,3,3).\]

$\blacksquare$ The case where one of the groups $\mathbb{W}_{i}$
($i\in\{1,2\}$) of Lemma \ref{lemme produit de 2 refl} is equal
to $G(3,3,3)$ has already been studied. In that case\[
\mathbb{W}_{1}\times\mathbb{W}_{2}\simeq G(3,3,3)\times G(3,3,2)\]
 is a reflection sub-group of $G_{S}$.

$\blacksquare$ If $[\,]^{2}\times[\,]^{2}$ is a reflection sub-group
of $G_{S}$, there exist coordinates such that : \[
F=\{x_{1}^{2}l_{1}(x_{3},x_{4},x_{5})+x_{2}^{2}l_{2}(x_{3},x_{4},x_{5})+x_{3}^{3}+x_{4}^{3}+x_{5}^{3}-3\lambda x_{3}x_{4}x_{5}=0\}\]
where $l_{1}$ and $l_{2}$ are two linearly independent forms. The
points $\mathbb{C}e_{1}$ and $\mathbb{C}e_{2}$ are vertices of cones
on $F$. The Fano surface $S$ contains two elliptic curves $E,E'$
such that $EE'=1$. 

$\blacksquare$ If $G(3,3,2)\times[\,]^{2}$ is a reflection sub-group
of $G_{S}$, then there exist coordinates such that:\[
F=\{x_{1}^{3}+x_{2}^{3}-3\lambda x_{1}x_{2}l_{1}(x_{4},x_{5})+x_{3}^{2}l_{2}(x_{4},x_{5})+x_{4}^{3}+x_{5}^{3}=0\},\]
where $l_{1}$ are $l_{2}$ linear forms. The points $\mathbb{C}(e_{1}-\beta e_{2}),\beta\in\mu_{3}$
and $\mathbb{C}e_{3}$ are vertices of cones. The Fano surface contains
three disjoint elliptic curves that cut another elliptic curve. The
dual graph of this configuration is the graph $D_{4}$.

$\blacksquare$ The last case is the group $G(3,3,2)\times G(3,3,2)$
for which the cubic is: \[
F=\{x_{1}^{3}+x_{2}^{3}+3ax_{1}x_{2}x_{5}+x_{3}^{3}+x_{4}^{3}+3bx_{3}x_{4}x_{5}+x_{5}^{3}=0\},\]
$(a,b\in\mathbb{C}$). The points $\mathbb{C}(e_{1}-\beta e_{2}),\,\mathbb{C}(e_{3}-\beta e_{4}),\,\beta\in\mu_{3}$
are vertices of cones on $F$. The Fano surface $S$ contains $3$
disjoint elliptic curves $E_{12}^{\beta},\,\beta\in\mu_{3}$ isomorphic
to the plane cubic \[
\{x_{3}^{3}+x_{4}^{3}+(1+a^{3})x_{5}^{3}+3bx_{3}x_{4}x_{5}=0\},\]
 and three others disjoint elliptic curves $E_{34}^{\alpha},\,\alpha\in\mu_{3}$
isomorphic to the cubic:\[
\{x_{1}^{3}+x_{2}^{3}+(1+b^{3})x_{5}^{3}+3ax_{1}x_{2}x_{5}=0\}.\]
These curves verify $E_{12}^{\beta}E_{34}^{\gamma}=1$ ($\beta,\gamma\in\mu_{3}$).
The surface $S$ contains $6$ elliptic curves.

Recall that for each elliptic curve $E\hookrightarrow S$, we denoted
by $\sigma_{E}:S\rightarrow S$ the involution associated to $E$
and that $G_{S}$ is a reflection group that is isomorphic to the
group ${\bf G}_{S}$ generated by the automorphisms $\sigma_{E}$.
Let us denote by $n_{s}$ the number of elliptic curves on $S$. The
following theorem gives the configuration classification of elliptic
curves on Fano surfaces:
\begin{thm}
\label{THE CLASSIFICATION THEOREM}(Classification Theorem). If the
group $G_{S}$ is irreducible, then it is isomorphic to one of the
following groups:\\
\begin{tabular}{|c|c|c|c|c|c|c|}
\hline 
Group $G_{S}$ & $\{1\}$ & $[\,]^{2}$ & $G(3,3,2)$ & $\Sigma_{4}$ & $\Sigma_{5}$ & $G(3,3,5)$\tabularnewline
\hline 
$n_{S}$ & 0 & 1 & 3 & 6 & 10 & 30\tabularnewline
\hline 
Order of $G_{S}$ & 1 & 2 & 6 & 24 & 120 & 9720\tabularnewline
\hline
\end{tabular}

Otherwise, $G_{S}$ is isomorphic to one of the following groups:\\
\begin{tabular}{|c|c|c|c|}
\hline 
$[\,]^{2}\times[\,]^{2}$ & $G(3,3,2)\times[\,]^{2}$ & $G(3,3,2)\times G(3,3,2)$ & $G(3,3,3)\times G(3,3,2)$\tabularnewline
\hline 
2 & 4 & 6 & 12\tabularnewline
\hline 
4 & 12 & 36 & 324\tabularnewline
\hline
\end{tabular}\\
Let $E,E'$ be two elliptic curves on $S$. We know the intersection
number $EE'$ and a plane model of $E$. 
\end{thm}

\section{Intermediate Jacobians isomorphic to a product of elliptic curves.}

Let $\lambda\in\mathbb{C}$, $\lambda^{3}\not=1$. The Fano surface
$S_{\lambda}$ of the cubic \[
F_{\lambda}=\{x_{1}^{3}+x_{2}^{3}+x_{3}^{3}-3\lambda x_{1}x_{2}x_{3}+x_{4}^{3}+x_{5}^{3}\}\hookrightarrow\mathbb{P}^{4}\]
 possesses $12$ smooth curves of genus $1$ for which we use the
notations of the previous paragraph. Let $\alpha$ be a primitive
third root of unity. Let $\lambda\in\mathbb{C}$ be such that the
elliptic curve $E_{\lambda}=\{x_{1}^{3}+x_{2}^{3}+x_{3}^{3}-3\lambda x_{1}x_{2}x_{3}=0\}$
has complex multiplication by $\mathbb{Q}(\alpha)$. Let $A_{\lambda}$
be the Albanese variety of the Fano surface $S_{\lambda}$ of $F_{\lambda}$
and let $\vartheta:S_{\lambda}\rightarrow A_{\lambda}$ be an Albanese
map. In this last section, we prove the following result:
\begin{thm}
\label{il existe une infinit=0000E9 de NS de rang 25}The Abelian
variety $A_{\lambda}$ is isomorphic to a product of elliptic curves. 
\end{thm}
Let $S$ be a Fano surface. The Albanese map $\vartheta:S\rightarrow A$
is an embedding and we consider points of $S$ as points of $A$.
\begin{lem}
\cite[(11.9)]{Clemens} There is a point $u_{o}$ on $A$ such that
for all points $s_{1},s_{2},s_{3}$ on $S$ such that the lines $L_{s_{1}}$,
$L_{s_{2}}$, $L_{s_{3}}$ are coplanar, we have:\[
s_{1}+s_{2}+s_{3}=u_{o}.\]

\end{lem}
Let $E$ be an elliptic curve on a Fano surface $S$. In order to
avoid heavy notations, we consider points of $E$ as points of $A$
and the morphism of Abelian varieties are taken modulo translation
: the differential and degree of such class is well defined.\\
For the points $s$ on $S$, the lines $L_{s}$, $L_{\gamma_{E}s}$,
$L_{\sigma_{E}s}$ are coplanar, hence the morphism $s\rightarrow s+\sigma_{E}s+\gamma_{E}s$
is constant i.e. its differential is $0$. Let $I_{A}$ be the identity
of $A$, let \[
\Gamma_{E}:A\rightarrow E,\,\Sigma{}_{E}:A\rightarrow A\]
 be the morphisms such that $\Gamma_{E}\circ\vartheta=\gamma_{E}$
and $\Sigma{}_{E}\circ\vartheta=\vartheta\circ\sigma_{E}$. The differential
of $\Sigma_{E}$ is $\sigma_{E}^{*}$. The morphism $Id_{A}+\Sigma_{E}+\Gamma_{E}$
is constant, hence: 
\begin{lem}
\label{lemme differentielle de la fibration}The differential of $\Gamma_{E}$
is the endomorphism $N_{E}\in End(H^{o}(\Omega_{S})^{*})$ such that
: \[
Id+\sigma_{E}^{*}+N_{E}=0\]
where $I$ is the identity of $H^{o}(\Omega_{S})^{*}$.
\end{lem}
Let us take $S=S_{\lambda}$. We denote by $\Gamma_{ij}^{\beta}:A_{\lambda}\rightarrow E_{ij}^{\beta}$
the morphism such that $\Gamma_{ij}^{\beta}\circ\vartheta=\gamma_{E_{ij}^{\beta}}$. 
\begin{lem}
\label{A  est isog=0000E8ne =0000E0 un produit de c ell}The degree
of the morphism \[
\Gamma=(\Gamma_{12}^{1},\Gamma_{23}^{1},\Gamma_{12}^{\alpha},\Gamma_{45}^{1},\Gamma_{45}^{\alpha}):A_{\lambda}\rightarrow E_{12}^{1}\times E_{23}^{1}\times E_{12}^{\alpha}\times E_{45}^{1}\times E_{45}^{\alpha}\]
 divides $81$.\end{lem}
\begin{proof}
Let $\Upsilon$ be the morphism $\Gamma$ composed with the natural
morphism $E_{12}^{1}\times E_{23}^{1}\times E_{12}^{\alpha}\times E_{45}^{1}\times E_{45}^{\alpha}\rightarrow A_{\lambda}$.
By Lemma \ref{lemme differentielle de la fibration}, we can compute
the differential $d\Upsilon$ of $\Upsilon$ and we find $|\det(d\Upsilon)|^{2}=81$,
hence the assertion.
\end{proof}
We can now complete the proof of Proposition \ref{il existe une infinit=0000E9 de NS de rang 25}:
\begin{proof}
By Lemma \ref{A  est isog=0000E8ne =0000E0 un produit de c ell} the
Abelian variety $A_{\lambda}$ is isogenous to a product of elliptic
curves ; by the choice of the parameter $\lambda$, these curves have
complex multiplication by the same field : the Néron-Severi group
of $A_{\lambda}$ has thus rank $25=h^{1,1}(A_{\lambda})$. Then \cite[Chap. 5, Exer. 5.6 (10)]{Birkenhake}
imply that $A_{\lambda}$ is isomorphic to a product of elliptic curves.
\end{proof}

Xavier Roulleau~~~~~roulleau@ms.u-tokyo.ac.jp\\
Graduate school of mathematical sciences, The University of Tokyo
3-8-1 Komaba, Meguro, Tokyo, 153-8914 Japan 
\end{document}